\documentclass[conference, 10 pt, conference]{ieeeconf}  % Comment this line out
\IEEEoverridecommandlockouts                              % This command is only
\usepackage{float}
\usepackage{amsmath}

\usepackage{amssymb}
\usepackage{algorithm} 
\usepackage{algpseudocode} 
\usepackage{xcolor}
\usepackage{bbm}     
\usepackage{graphicx} 
\usepackage[normalem]{ulem}
\newtheorem{theorem}{Theorem}
\newtheorem{lemma}{Lemma}

\newtheorem{definition}{Definition}
\newtheorem{problem}{Problem}
\usepackage{dsfont}

\title{\LARGE \bf Data-Driven Stable Neural Feedback Loop Design

}

\author{Zuxun Xiong, Han Wang, Liqun Zhao, and Antonis Papachristodoulou% <-this % stops a space
\thanks{Z. Xiong, H. Wang, Liqun Zhao, and A. Papachristodoulou are with Department of Engineering Science, University of Oxford, Parks Road, Oxford, OX1 3PJ, U.K. {\tt\small \{zuxun.xiong, han.wang, liqun.zhao,antonis\}@eng.ox.ac.uk }}
\thanks{AP was supported in part by UK's Engineering, Physical Sciences Research Council projects EP/X017982/1 and EP/Y014073/1.}
}

\begin{document}

\maketitle
\thispagestyle{empty}
\pagestyle{empty}

%%%%%%%%%%%%%%%%%%%%%%%%%%%%%%%%%%%%%%%%%%%%%%%%%%%%%%%%%%%%%%%%%%%%%%%%%%%%%%%%
\begin{abstract}
This paper proposes a data-driven approach to design a feedforward Neural Network (NN) controller with a stability guarantee for plants with unknown dynamics. We first introduce data-driven representations of stability conditions for Neural Feedback Loops (NFLs) with linear plants, which can be formulated into a semidefinite program (SDP). Subsequently, this SDP constraint is integrated into the NN training process to ensure stability of the feedback loop. The whole NN controller design problem can be solved by an iterative algorithm. Finally, we illustrate the effectiveness of the proposed method compared to model-based methods via numerical examples. 

\end{abstract}

%%%%%%%%%%%%%%%%%%%%%%%%%%%%%%%%%%%%%%%%%%%%%%%%%%%%%%%%%%%%%%%%%%%%%%%%%%%%%%%%
\section{Introduction}
Decades before the recent rapid development of Artificial Intelligence (AI), pioneering research had already proposed the application of neural networks (NNs) as controllers for systems with unknown dynamics \cite{RN102, RN101, RN72}. As related technologies have advanced, this research topic is now relevant in feedback control and has achieved success in several applications in recent years, e.g., 
% \cite{RN102} and \cite{RN101} used a single-layer NN and a double-layer NN to stablilize an inverted pendulum respectively. NN controller based on measurement data obtained during the system operation \cite{RN72}. 
optimal control \cite{RN75}, model predictive control \cite{RN76} and reinforcement learning \cite{RN13}. 

The ability of NNs to act as general function approximators underpins their impressive performance on control. However, this characteristic also presents challenges in providing rigorous stability and robustness guarantees. To address this issue, research in NN verification focuses on the relationship between inputs and outputs of NNs, ensuring they can operate reliably. In \cite{RN23}, an interval bound propagation (IBP) approach was proposed to calculate an estimate on the output range of every NN layer. This method is relatively simple but also very conservative. Quadratic constraints (QC) were later used to bound the nonlinear activation functions in NNs \cite{RN40}: the NN verification problem with these QC bounds was formulated as a semidefinite program and solved efficiently. In \cite{RN70}, a tighter bound consisting of two sector constraints for different activation functions with higher accuracy was proposed. A Lipschitz bound estimation method \cite{RN39} can also be used to evaluate the robustness of NNs. 

Based on these bounds and stability conditions provided by the NN verification literature, recent results focused on the design problem of NN feedback control. By imposing a Lipschitz bound constraint in the training process, a robust NN controller was designed \cite{RN27}. In \cite{RN15}, a framework to design a stable NN controller through imitation learning (IL) for linear time-invariant (LTI) systems was proposed. To ensure that controllers can retain and process long-term memories, recurrent neural network controllers were trained with stability guarantees for partially observed linear systems \cite{RN54} and nonlinear systems \cite{RN28}.
% Although the imitation learning framework of \cite{RN15} can train NN controller with stability guarantee, it
However, most of the existing studies assume the dynamics of the controlled systems are known, or at least partially known. This assumption contradicts the original intention of designing NN controllers, which is to control systems with unknown dynamics. This is exactly the gap we want to fill in this work. 

% With technological advances on computational and storage capacity, as well as more access to data, data-driven control for complex systems can facilitate research in both control and other research fields \cite{RN12}. 
Data-driven control method promises to address this problem \cite{RN12}. A thorough introduction on how to apply data-driven methods to represent linear systems and design stable feedback control systems was given in \cite{RN10}. Relevant methods were applied to analyse the stability of nonlinear systems \cite{RN71} using  Sum of Squares (SOS) \cite{RN87}. We also proposed a data-driven verification method to analyze the closed-loop system with an NN controller \cite{wang2023model}. In this study, we use the term Neural Feedback Loop (NFL) to represent the closed-loop feedback system with an NN controller.

The first contribution of this paper is a data-driven representation of stability conditions for an NFL with unknown linear plant. Based on the stability conditions, we propose an iterative design algorithm to train a stable NN. Secondly, an NN fine-tuning algorithm is proposed to stabilize an existing NN controller for an unknown system. Both these algorithms require only the collection of data generated by persistently exciting inputs.

% Finally, in case studies, we train NN controllers through imitation learning based on the proposed design algorithm. We also fine-tune an NN controller trained by imitation learning using our fine-tuning algorithm. The effectiveness of both algorithms are verified. 

The rest of this paper is organized as follows. In Section \ref{sec:pre}, we review tools for NFL analysis and the data-driven representation of linear state-feedback systems. In Section \ref{sec:Analysis}, we formulate a data-driven stability condition for NFLs. Based on these conditions, we propose a stable NFL design algorithm and an NN fine-tuning algorithm in Section \ref{sec:Algorithm}, and then illustrate two algorithms using numerical examples in Section \ref{sec:case}. We draw conclusions in the last section. 

\textbf{Notation:} We use $\mathbb{R}^{n,m}$,  $\mathbb{S}^{n}_{++}$ and $I_n$ to denote $n$-by-$m$ dimensional real matrices, $n$-by-$n$ positive definite matrices and $n$-by-$n$ identity matrices, respectively. We use $||\cdot||_F$ and tr$(\cdot)$ to denote the Frobenius norm and the trace of a matrix.
 
% \subsection{Sum of Square and Stability Analysis}
% We will refer to \cite{RN69} and \cite{RN61}.

%%%%%%%%%%%%%%%%%%%%%%%%%%%%%%%%%%%%%%%%%%%%%%%%%%%%%%%%%%%%%%%%%%%%%%%%%%%%%%%%

\section{Preliminaries}\label{sec:pre}
% {\color{blue}This is missing, some notations, e.g. $I_n$, $1_n$, $\mathbb{S}_{++}^n$, are used without definitions.}
% This section provides background information used in the subsequent sections. First of all, we give a basic introduction of NFL. Then we illustrate how to verify the stability of NFL based on the works mentioned before. Finally, fundamental formulas for data-driven control method are briefly introduced. 
% we briefly introduce the formulas for data-driven control method. Finally, several conclusions on NFL stability verification from our previous work \cite{wang2023model} are illustrated. 
%Finally, the sum of square (SOS) programming, which will provide a foundation of stability analysis on nonlinear system, will be introduced. 
\subsection{Neural Feedback Loop}
A generic NFL is a closed-loop dynamical system consisting of a plant $G$ and an NN controller $\pi(\cdot)\in\mathbb{R}^{n_\pi}$. In this work, we consider $G$ to be an LTI system of the form:
\begin{equation}\label{eq:system}
    x(k+1) = A_Gx(k)+B_Gu(k),
\end{equation}
where $A_G\in\mathbb{R}^{n_x \times n_x}$ and $B_G\in\mathbb{R}^{n_x\times n_\pi}$. Here $x(k)\in\mathbb{R}^{n_x}$ and $u(k)\in\mathbb{R}^{n_{\pi}}$ are the system state and the control input, respectively. The controller is designed as a feed-forward fully connected neural network $\pi(k)$ with $l$ layers as follows:\begin{subequations}\label{eq:nn}
    \begin{eqnarray}
        \omega^0(k) &=& x(k),\\
        \nu^i(k) &=& W^i\omega^{i-1}(k)+b^i,~\mathrm{for}~i=1,\ldots,l,\\
        \omega^i(k)&=&\phi^i(\nu^{i}(k)),~\mathrm{for}~i=1,\ldots,l,\\
        \pi(k)&=&W^{l+1}\omega^l(k)+b^{l+1}.
    \end{eqnarray}
\end{subequations}
For the $i^{\mathrm{th}}$ layer of the NN, we use $W^i\in\mathbb{R}^{n_{i}\times n_{i-1}}$,  $b^{i}\in\mathbb{R}^{n_i}$, $n_i$, $\nu^i$ and $\omega^i$ to denote its weight matrix, bias vector, number of neurons, and corresponding vectorized input and output, respectively. The input of the NN controller is $\omega^0(k)=x(k)$, which is the state of the plant at time $k$. $\phi^i(\cdot):\mathbb{R}^{\nu_i}\to \mathbb{R}^{\nu_i}$ is a vector of nonlinear activation functions on the $i^{\mathrm{th}}$ layer, defined as
\begin{equation}
    \phi^i(\nu^i) = [\varphi(\nu^i_1),\ldots,\varphi(\nu^i_{n_i})]^\top,
\end{equation}
where $\varphi(\cdot):\mathbb{R}\to \mathbb{R}$ is an activation function, such as ReLU, sigmoid, and $\tanh$.
% such as ReLU $\varphi(\nu) := \max(0,\nu)$, tanh $\varphi(\nu):=\tanh(\nu)$, sigmoid $\varphi(\nu) :=\frac{1}{1+e^{-\nu}}$. 

\subsection{Stability Verification for Neural Feedback Loop}
To verify stability of an NFL, a commonly used method is to isolate the nonlinear terms introduced by the activation functions, then treat the nonlinearity as a disturbance \cite{RN15}. The NFL in \eqref{eq:nn} can be rewritten as:
\begin{subequations}\label{eq:NN-linear}
    \begin{align}
        \begin{bmatrix}
            \pi(k)\\
            \nu_\phi(k)
        \end{bmatrix}&=
         N \begin{bmatrix}
            x(k)\\
            \omega_\phi(k)
        \end{bmatrix}\\
        \omega_\phi(k)&=\phi(\nu_\phi(k)).\label{eq:nn-nonlinear}
    \end{align}
\end{subequations}
where $\nu_{\phi}(k)\in \mathbb{R}^{n_\phi}$ and $\omega_{\phi}(k) \in \mathbb{R}^{n_\phi}$ are vectors formed by stacking $\nu^i(k)$ and $\omega^i(k)$ of each layer $i$, $i=1,\ldots,l$, $\phi(\cdot): \mathbb{R}^{n_\phi} \rightarrow \mathbb{R}^{n_\phi}$ is the stacked activation function for all layers, where $n_\phi:=\sum^{l}_{i=1}n_i$. $N:=
    \begin{bmatrix}
        N_{\pi x}&N_{\pi\omega}\\
        N_{\nu x}&N_{\nu \omega}
    \end{bmatrix}$ is a matrix consisting of NN weights $W^i$.  
%sector constraints
It should be noticed that we impose the constraint $\pi(0)=0$ to ensure the equilibrium remains at the origin with the NN controller. To achieve this, we also set all bias $b^i$ to be 0. 
To deal with the nonlinearities in \eqref{eq:nn-nonlinear}, sector constraints have been proposed to provide lower and upper bounds for different kinds of activation functions of NNs. The interested readers are referred to \cite{RN40,RN70} for a comprehensive review and comparison on different sectors. In this work, we rely on a commonly used local sector:
\begin{definition}
 Let $\alpha, \beta, \underline{\nu}, \bar{\nu} \in\mathbb{R}$ with $\alpha \leq \beta$ and $\underline{\nu} \leq 0 \leq \bar{\nu}$. The function $\varphi:\mathbb{R}\to \mathbb{R}$ satisfies the \emph{local sector} $[\alpha,\beta]$ if 
    \begin{equation}
        ( \varphi(\nu)-\alpha \nu) \cdot(\beta \nu- \varphi(\nu)) \geq 0, \quad \forall \nu \in[\underline{\nu}, \bar{\nu}].
    \end{equation}
\end{definition}
For example, $\varphi(\nu):=\tanh(\nu)$ restricted to the interval [$-\bar{\nu},\bar{\nu}$] satisfies the local sector $[\alpha,\beta]$ with $\alpha=\tanh(\bar{\nu})/\bar{\nu}>0$ and $\beta=1$.

% \begin{equation}\notag
%     \alpha:=\text{min} (\frac{\text{tanh}(\bar{\nu})-\text{tanh}(\nu_*)}{\bar{\nu}-\nu_*},  \frac{\text{tanh}(\nu_*)-\text{tanh}(\underline{\nu})}{\nu_*-\underline{\nu}}).
% \end{equation}

By composing the sector bounds for each individual activation function, we can obtain sector constraints for the stacked nonlinearity $\phi(\cdot)$. One composition method is shown in the following lemma.
\begin{lemma}[\protect {\cite[Lemma 1]{RN15}}]\label{lem:sector}
    Let $\alpha_\phi$, $\beta_{\phi}, \underline{\nu},  \bar{\nu} \in\mathbb{R}^{n_\phi}$ be given with $\alpha_{\phi}\le \beta_{\phi}$, and $\underline{\nu} \leq 0 \leq \bar{\nu}$. Assume $\phi$ element-wisely satisfies the local sector $[\alpha_{\phi},\beta_{\phi}]$ for all $\nu_\phi\in [\underline{\nu}, \bar{\nu}]$. Then, for any $\lambda\in\mathbb{R}^{n_\phi}$ with $\lambda\ge 0$, and for all $\nu_\phi\in [\underline{\nu}, \bar{\nu}]$, $\omega_\phi=\phi(\nu_\phi)$, we have
    \begin{align}\label{eq:sector}
        \begin{bmatrix}
            \nu_\phi\\
            \omega_\phi
        \end{bmatrix}^\top
        \begin{bmatrix}
            -2A_\phi B_\phi\Lambda&(A_\phi+B_\phi)\Lambda\\
            (A_\phi+B_\phi)\Lambda&-2\Lambda
        \end{bmatrix}
        \begin{bmatrix}
            \nu_\phi\\
            \omega_\phi
        \end{bmatrix}\ge 0,
    \end{align}
    where $A_\phi$ = diag($\alpha_\phi$), $B_\phi$ = diag($\beta_\phi$), and $\Lambda$ = diag($\lambda$). 
\end{lemma}
%stability constraints
% The computational method of $A_{\phi}$ and $B_{\phi}$, i.e., the sector bound of every layer of NN, can be found in \cite{RN23}.
Putting the result from Lemma \ref{lem:sector} into a robust Lyapunov stability analysis framework, we can obtain a sufficient condition for an NFL to be stable.
\begin{theorem}[\protect {\cite[Theorem 1]{RN15}}]\label{theo:stability}
Consider an NFL with plant $G$ satisfying \eqref{eq:system} and NN controller $\pi$ as in \eqref{eq:nn} with an equilibrium point $x_*=0_{n_x}$ and a state constraint set $X \subseteq\{x: -\bar{x} \leq Hx \leq \bar{x}\}$. 
% Let $\underline{\nu}, \bar{\nu}$ be corresponding activation input bounds for state $x$ and 
% Assume $\phi(\cdot)$ element-wisely satisfies the offset sector $[\alpha_\phi,\beta_\phi]$, centered on the point $(\nu_*,\omega_*)$, for all $\nu_\phi \in [\underline{\nu},\bar{\nu}]$.
If there exist a positive definite matrix $P\in\mathbb{S}^{n_x}_{++}$, a vector $\lambda \in \mathbb{R}^{n_\phi}$ with $\lambda\geq0$, and a matrix $\Lambda := \text{diag}(\lambda)$ that satisfy
\begin{subequations}\label{eq:original_stability}
\begin{equation}\label{eq:original_stability1}
\begin{aligned}
R_V^{\top} & {\left[\begin{array}{cc}
A_G^{\top} P A_G-P & A_G^{\top} P B_G \\
B_G^{\top} P A_G & B_G^{\top} P B_G
\end{array}\right] R_V } \\
& +R_\phi^{\top}\left[\begin{array}{cc}
-2 A_\phi B_\phi \Lambda & \left(A_\phi+B_\phi\right) \Lambda \\
\left(A_\phi+B_\phi\right) \Lambda & -2 \Lambda
\end{array}\right] R_\phi \prec  0,
\end{aligned}
\end{equation}
and
\begin{equation}\label{eq:state constraint}
    \begin{bmatrix}
            \bar{x}^2_i & H^\top_i\\
            H_i & P
        \end{bmatrix} \geq 0, i=1,\ldots,n_x,
\end{equation}
where 
\begin{equation}\label{eq:RV RPhi}
    R_V:=\left[\begin{array}{cc}
I_{n_x} & 0_{n_x \times n_\phi} \\
N_{\pi x} & N_{\pi \omega}
\end{array}\right], R_\phi:=\left[\begin{array}{cc}
N_{\nu x} & N_{\nu \omega} \\
0_{n_\phi \times n_x} & I_{n_\phi}
\end{array}\right],
\end{equation}
\end{subequations}
where $H_i^\top$ is the $i^{th}$ row of the matrix $H$, then the NFL is locally asymptotically stable around the equilibrium point $x_*$, and the ellipsoid $\mathcal{E}(P):=\left\{x \in \mathbb{R}^{n_x}: x^{\top} P x \leq 1\right\}$ is an inner-approximation to the region of attraction (ROA).
\end{theorem}

% The first term and second term of the left side of condition \eqref{eq:original_stability1} come from Lyapunov function $V(x):=(x-x_*)^{\top}P(x-x_*)$ and sector constraint \eqref{eq:sector} respectively.

\subsection{Data-Driven Representation of the System}
Verifying stability of an NFL using \eqref{eq:original_stability} requires $A_G$ and $B_G$ to be precisely known. We refer to this type of approaches \emph{model-based}. In contrast to model-based approaches, direct \emph{data-driven} approaches aim to bypass the identification step and use data directly to represent the system. Consider System \eqref{eq:system}. We carry out experiments to collect $T$-long time series of system inputs, states, and successor states as follows:
\begin{subequations}
    \begin{align}
    &U_{0,T}:=\begin{bmatrix}
        u(0)&u(1)&\ldots&u(T-1)
    \end{bmatrix}\in\mathbb{R}^{n_u\times T}\\
    &X_{0,T}:=\begin{bmatrix}
        x(0)&x(1)&\ldots&x(T-1)
    \end{bmatrix}\in\mathbb{R}^{n_x\times T}\\
    &X_{1,T}:=\begin{bmatrix}
        x(1)&x(2)&\ldots&x(T)
    \end{bmatrix}\in\mathbb{R}^{n_x\times T}
\end{align}
\end{subequations}
It should be noted that here we use $u(t)$ to denote the instantaneous input signal for the open-loop system. The signal is not necessary produced by the NN controller $\pi(\cdot)$. 
% We assume that the data are collected in an accurate way without noise corruptions. 
This data series can be used to represent any $T$-long trajectory of the linear dynamical system as long as the following rank condition holds \cite{RN10}:
\begin{equation}\label{eq:rank-condition}
    \mathrm{rank}\left(\begin{bmatrix}
        U_{0,T}\\X_{0,T}
    \end{bmatrix}\right)=n_u+n_x.
\end{equation}
To satisfy this rank condition, the collected data series should be ``sufficiently long", i.e., $T \geq (n_u+1)n_x+n_u$. Under the rank condition, System \eqref{eq:system} with a state-feedback controller $u(k)=Kx(k)$ has the following data-driven representation \cite[Theorem 2]{RN10}:
\begin{equation}
    x(k+1)=(A_G+B_GK)x(k)=X_{1,T}G_Kx(k),
\end{equation}
where $G_K$ is a $T \times n_x$ matrix satisfying
\begin{equation}
    \begin{bmatrix}
        K \\
        I_{n_x} \\
    \end{bmatrix}
    =
    \begin{bmatrix}
        U_{0,T} \\
        X_{0,T} \\
    \end{bmatrix}G_K.
\end{equation}
% and therefore 
% \begin{equation}
%     u(k)=U_{0,T}G_Kx(k).
% \end{equation}

In our problem, the controller $\pi(\cdot)$ is a highly nonlinear NN. This makes it challenging to design a stable NFL directly from data.
We will demonstrate how to combine this method with the NFL stability condition in the sequel.

%%%%%%%%%%%%%%%%%%%%%%%%%%%%%%%%%%%%%%%%%%%%%%%%%%%%%%%%%%%%%%%%%%%%%%%%%%%%%%%%
% \section{Problem Formulation}\label{sec:problem}
\subsection{Problem Formulation}
Our problems of interest are:
\begin{problem}\label{prob:1}
    Consider an unknown plant $G$ as in \eqref{eq:system}. Design a feed-forward fully connected NN controller $\pi:\mathbb{R}^{n_x}\to \mathbb{R}^{n_\pi}$ that optimises a given loss function and stabilises the plant around the origin. 
\end{problem}

\begin{problem}\label{prob:2}
    Consider an unknown plant $G$ as in \eqref{eq:system} and a \emph{given} NN controller $\pi:\mathbb{R}^{n_x}\to \mathbb{R}^{n_\pi}$. Minimally tune the NN to guarantee stability of the NFL around the origin.
\end{problem}

% {\color{blue}HW: I suggest to simplify and remove the following discussions. They are a bit crowded at this stage. But I'm happy with keep them. Let's discuss.}
% \sout{Consider an NFL consisting of an unknown plant $G$ as \eqref{eq:system} and a feed-forward fully connected NN controller $\pi:\mathbb{R}^{n_x}\to \mathbb{R}^{n_\pi}$ as \eqref{eq:nn}, which satisfies $\pi(0)=0$. How to design an NN controller, that is, determine the values of $W^i$ and $b^i$, so that the system can be stabilized while some performance metrics are optimized, e.g., maximize the size of region of attraction (ROA) or minimize the loss of the NN? Furthermore, due to the difficulty in obtaining parameters of the unknown system $G$ through traditional approach, the design method of the stable NFL needs to be data-driven.}

% \sout{Besides designing new NN controllers for unknown systems, we also note a class of application scenarios suitable for the proposed method: that is, by fine-tuning the parameters of a given NN controller, making an originally unstable system stable. Compared to introducing a new safety filter to the NFL \cite{RN97}, fine-tuning the parameters of the existing NN might be more convenient. So, how to fine-tuning existing NN becomes our second problem. }

Throughout this paper, $\tanh$ is used as the activation function $\varphi(\cdot)$ of the NN $\pi(\cdot)$. Our result could be easily extended to other types of activation functions, using different sectors.
% \addtolength{\textheight}{-3cm}   % This command serves to balance the column lengths

%%%%%%%%%%%%%%%%%%%%%%%%%%%%%%%%%%%%%%%%%%%%%%%%%%%%%%%%%%%%%%%%%%%%%%%%%%%%%%%%
 
\section{Data-Driven Stability Analysis for NFLs}\label{sec:Analysis}
In this section we provide data-driven stability conditions for NFLs. To obtain convex conditions, a loop transformation is used to normalise the sector constraints~\cite{RN15}.% first give a review on the stability conditions on both linear system and nonlinear system. Then we give our theorems and proofs on data-driven representation of the stability conditions for linear NFL and nonlinear NFL respectively. 

\subsection{Loop Transformation}
With a given controller $\pi$, we can calculate the bounds $\underline{\nu}, \bar{\nu}$ for every layer's output based on a given state constraint set. If the plant $G$ is also given, then  $A_G$, $B_G$, and $N$ are known and the stability condition \eqref{eq:original_stability} is convex in matrices $P$ and $\Lambda$. Stability can then be \emph{verified} by solving an SDP. However, for the design problem \ref{prob:1}, $G$ is unknown and $\pi$ is to be designed, which means that $A_G$, $B_G$ and $N$ are all decision variables in the problem \eqref{eq:original_stability}. The condition is then nonconvex. To deal with the nonlinearity, we first utilize a loop transformation to normalise the sector constraints.
\begin{subequations}\label{eq:NN-loop-transform}
    \begin{align}
        \begin{bmatrix}
            \pi(k)\\
            \nu_\phi(k)
        \end{bmatrix}&=
        \tilde N \begin{bmatrix}
            x(k)\\
            z_\phi(k)
        \end{bmatrix}\\
        z_\phi(k)&=\tilde \phi(\nu_\phi(k)).
    \end{align}
\end{subequations}
The new internal state $z_\phi(k)$ is related to $\omega_\phi(k)$, as follows
\begin{equation}
    \omega_\phi(k)=\frac{B_\phi-A_\phi}{2} z_\phi(k)+\frac{A_\phi+B_\phi}{2} \nu_\phi(k).
\end{equation}
Through the transformation, the nonlinearity $\tilde \phi$ is normalised, namely, lies in sector $[-1_{n_\phi \times 1}, 1_{n_\phi \times 1}]$. Using Lemma \ref{lem:sector}, we have
\begin{equation}
    \begin{bmatrix}
        \nu_\phi\\z_\phi
    \end{bmatrix}^\top 
    \begin{bmatrix}
        \Lambda&0\\0&-\Lambda
    \end{bmatrix}
    \begin{bmatrix}
        \nu_\phi\\z_\phi
    \end{bmatrix}\ge 0,\quad \forall \nu_\phi\in[\underline{\nu},\overline{\nu}].
\end{equation}
The transformed matrix $\tilde N$ can be derived by
\begin{equation}\label{eq:f(N)compute}
    \tilde N  =\begin{bmatrix}
        N_{\pi x}+C_2(I-C_4)^{-1}N_{\nu x}&C_1+C_2(I-C_4)^{-1}C_3
        \\(I-C_4)^{-1}N_{\nu x}&(I-C_4)^{-1}C_3
    \end{bmatrix}
\end{equation}
% \begin{equation}\label{eq:f(N)compute}
% \begin{aligned}
%     \tilde{N} & =
%     \left[\begin{array}{ccc}
%         N_{\pi x}+C_2(I-C_4)^{-1}N_{\nu x} & C_1+C_2(I-C_4)^{-1}C_3 \\
%         (I-C_4)^{-1}N_{\nu x} & (I-C_4)^{-1}C_3 \\
%     \end{array}\right.\\
%     &\quad \left.\begin{array}{c}
%         N_{\pi b}+C_2(I-C_4)^{-1}N_{\nu b} \\
%         (I-C_4)^{-1}N_{\nu b} \\
%     \end{array}\right]:=
%     \begin{bmatrix}
%         \tilde N_{\pi x} & \tilde N_{\pi z} & \tilde N_{\pi b} \\
%         \tilde N_{\nu x} & \tilde N_{\nu z} & \tilde N_{\nu b}
%     \end{bmatrix}
% \end{aligned}
% \end{equation}
where
\begin{equation}\label{eq:Ccompute}
\begin{aligned}
    &C_1=N_{\pi \omega}\frac{B_\phi-A_\phi}{2},C_2=N_{\pi \omega}\frac{A_\phi+B_\phi}{2},\\
    & C_3=N_{\nu \omega}\frac{B_\phi-A_\phi}{2},C_4=N_{\nu \omega}\frac{A_\phi+B_\phi}{2}.
\end{aligned}
\end{equation}
We define $\tilde N:=
    \begin{bmatrix}
        \tilde N_{\pi x}&\tilde N_{\pi z}   \\
        \tilde N_{\nu x}&\tilde N_{\nu z}
    \end{bmatrix}$. Using the new system representation, the stability condition \eqref{eq:original_stability1} and \eqref{eq:RV RPhi} can be reformulated as:
\begin{subequations}\label{eq:transformed_stability}
\begin{equation}\label{eq:transformed_stability1}
\begin{aligned}
    \tilde{R}_V^{\top} & \left[\begin{array}{cc}
A_G^{\top} P A_G-P & A_G^{\top} P B_G \\
B_G^{\top} P A_G & B_G^{\top} P B_G
\end{array}\right] \tilde{R}_V\\
&+\tilde{R}_\phi^{\top}\left[\begin{array}{cc}
\Lambda & 0 \\
0 & -\Lambda
\end{array}\right] \tilde{R}_\phi \prec  0,
\end{aligned}
\end{equation}
where 
\begin{equation}\label{eq:transformed_stability2}
    \tilde{R}_V:=\left[\begin{array}{cc}
I_{n_x} & 0 \\
\tilde{N}_{\pi x} & \tilde{N}_{\pi z}
\end{array}\right], \tilde{R}_\phi:=\left[\begin{array}{cc}
\tilde{N}_{\nu x} & \tilde{N}_{\nu z} \\
0 & I_{n_\phi}
\end{array}\right].
\end{equation}
\end{subequations}

\subsection{Data-Driven Representation of Stability Condition for NFLs}
% For simplicity, we will proceed with the theoretical derivation directly based on our previous work \cite{wang2023model}. Unlike the verification problem, in the design problem, the parameters of NN are no longer given but appear as decision variables in the problem. Consequently, the stability condition, (25a), in \textit{Theorem 4} of \cite{wang2023model}  will become:
The following theorem gives a sufficient condition for an NN controller $\pi$ to stabilize an NFL with an unknown plant $G$.

% {\color{blue}Is this theorem similar as the one in L4DC paper? If so, I don't think a proof is necessary here.}

% {\color{blue}The following Theorem should be revised.}
\begin{theorem}\label{theo:stability_DDC}
Consider an unknown LTI plant $G$ \eqref{eq:system} with an equilibrium point $x_*=0_{nx}$ and a state constraint set $X \subseteq\{x: -\bar{x} \leq Hx \leq \bar{x}\}$. Let rank condition \eqref{eq:rank-condition} hold. Find a matrix $Q_1\in\mathbb{S}^{n_x}_{++}$, a diagonal matrix $Q_2\in\mathbb{S}^{n_\phi}_{++}$, matrices $L_1\in\mathbb{R}^{T\times n_x}$, $L_2\in\mathbb{R}^{T\times n_\phi}$, $L_3\in\mathbb{R}^{n_\phi \times n_x}$ and $L_4\in\mathbb{R}^{n_\phi \times n_\phi}$, such that \eqref{eq:state constraint} is satisfied, and
\begin{subequations}\label{eq:stability-design}
\begin{equation}\label{eq:stability1-design}
     \begin{bmatrix}
        Q_1 & 0 & L_1^TX_{1,T}^T & L_3^T\\
        0 & Q_2 & L_2^TX_{1,T}^T & L_4^T\\
        X_{1,T}L_1 & X_{1,T}L_2 & Q_1 & 0\\
        L_3 & L_4 & 0 & Q_2\\
    \end{bmatrix} \succ 0.
\end{equation}
If we can design an NN to satisfy:
% \begin{equation}\label{eq:L_34}
%     L_3=\Tilde{N}_{\nu x}Q_1, \quad L_4=\Tilde{N}_{\nu z}Q_2,
% \end{equation}
% and
\begin{gather}
    \Tilde{N}Q=\bar{U}_{0,T}L, \label{eq:dd con1}\\
    \begin{bmatrix}
        I & 0\\
    \end{bmatrix}Q
    =
    \bar{X}_{0,T}L,
    \label{eq:dd con2}
\end{gather}
\end{subequations}
where $Q:=
\begin{bmatrix} Q_1 & 0 \\ 0 & Q_2 \\ \end{bmatrix}$, $L:=
\begin{bmatrix} L_1 & L_2 \\ L_3 & L_4 \\ \end{bmatrix}$, $\bar{U}_{0,T}:=
\begin{bmatrix} U_{0,T} & 0 \\ 0 & I \\ \end{bmatrix}$ and $\bar{X}_{0,T}:=
\begin{bmatrix} X_{0,T} & 0 \\ \end{bmatrix}$, then the designed NN controller can make the LTI system locally asymptotically stable around $x_*$. Besides, the set $\mathcal{E}(P)$ is an ROA inner-approximation for the system.
\end{theorem}
\begin{proof}
By substituting \eqref{eq:transformed_stability2} into \eqref{eq:transformed_stability1} and then utilising Schur complement, we can obtain stability condition:
\begin{equation}\label{eq:stability-Schur1}
    \begin{bmatrix}
        P&0&A_G^\top +\tilde N_{\pi x}^\top B_G^\top & \tilde N_{\nu x}^\top \\
        0&\Lambda&\tilde N_{\pi z}^\top B_G^\top & \tilde N_{\nu z}^\top \\
        A_G+B_G\tilde N_{\pi x} & B_G\tilde N_{\pi z} & P^{-1}&0\\
        \tilde N_{\nu x}&\tilde N_{\nu z}&0&\Lambda^{-1}
    \end{bmatrix}\succ 0.
\end{equation}
Here $A_G$, $B_G$ are unknown parameters which can be replaced using data-driven methods. Based on condition \eqref{eq:rank-condition}, we can find $G_1 \in\mathbb{R}^{T\times n_x} $ and $G_2 \in\mathbb{R}^{T\times n_\phi}$ that satisfy:
\begin{equation}\label{eq:data-original-con}
\begin{bmatrix}
\Tilde{N}_{\pi x} \\
I \\
\end{bmatrix}
=
\begin{bmatrix}
U_{0,T} \\
X_{0,T} \\
\end{bmatrix}G_1, \quad
\begin{bmatrix}
\Tilde{N}_{\pi z} \\
0 \\
\end{bmatrix}
=
\begin{bmatrix}
U_{0,T} \\
X_{0,T} \\
\end{bmatrix}G_2.
\end{equation}
Then we can formulate $A_G+B_G\tilde N_{\pi x}$ and $B_G\tilde N_{\pi z}$ as:
\begin{subequations}
\begin{equation} \label{eq:data rep1}
\begin{aligned}  
    & A_G+B_G\tilde N_{\pi x}=
    \begin{bmatrix}
        B_G & A_G\\
    \end{bmatrix}
    \begin{bmatrix}
        \Tilde{N}_{\pi x} \\
        I \\
    \end{bmatrix}\\
    &=
    \begin{bmatrix}
        B_G & A_G\\
    \end{bmatrix}
    \begin{bmatrix}
        U_{0,T} \\
        X_{0,T} \\
    \end{bmatrix}G_1
    =
    X_{1,T}G_1,
\end{aligned}
\end{equation} 
\begin{equation} \label{eq:data rep2}
    B_G\tilde N_{\pi z}=
    \begin{bmatrix}
        B_G & A_G\\
    \end{bmatrix}
    \begin{bmatrix}
        \Tilde{N}_{\pi z} \\
        0 \\
    \end{bmatrix}\\
    =
    X_{1,T}G_2.
\end{equation} 
\end{subequations}
Now we can replace the terms containing system parameters with these two terms. Besides, multiplying condition \eqref{eq:stability-Schur1} by $\text{diag}(P^{-1}, \Lambda^{-1}, I_{n_x}, I_{n_\phi})$ on both left and right, we obtain:
\begin{equation}\label{eq:stability-Schur2}
    \begin{bmatrix}
        P^{-1} & 0 & P^{-1}G_1^\top X_1^\top & P^{-1}\tilde N_{\nu x}^\top \\
        0 &\Lambda^{-1} & \Lambda^{-1}G_2^\top X_1^\top & \Lambda^{-1}\tilde N_{\nu z}^\top \\
        X_1G_1P^{-1} & X_1G_2\Lambda^{-1} & P^{-1} & 0\\
        \tilde N_{\nu x}P^{-1} & \tilde N_{\nu z}\Lambda^{-1} & 0 & \Lambda^{-1}
    \end{bmatrix}\succ 0.
\end{equation}
% \begin{equation}\label{eq:stability-Schur2}
%     \begin{bmatrix}
%         P^{-1} & 0 & P^{-1}(A_G^\top+\tilde N_{\pi x}^\top B_G^\top) & P^{-1}\tilde N_{\nu x}^\top \\
%         0 &\Lambda^{-1} & \Lambda^{-1}\tilde N_{\pi z}^\top B_G^\top & \Lambda^{-1}\tilde N_{\nu z}^\top \\
%         (A_G+B_G\tilde N_{\pi x})P^{-1} & B_G\tilde N_{\pi z}\Lambda^{-1} & P^{-1} & 0\\
%         \tilde N_{\nu x}P^{-1} & \tilde N_{\nu z}\Lambda^{-1} & 0 & \Lambda^{-1}
%     \end{bmatrix}\succ 0.
% \end{equation}
Let $Q_1=P^{-1}$, $Q_2=\Lambda^{-1}$, $L_1=G_1Q_1$, $L_2=G_2Q_2$, $L_3=\tilde{N}_{\nu x}Q_1$ and $L_4=\tilde{N}_{\nu z}Q_2$. Then stability condition \eqref{eq:stability-Schur2} can be expressed as \eqref{eq:stability1-design} and condition \eqref{eq:data-original-con} from data-driven method can be reformulated to \eqref{eq:dd con1} and \eqref{eq:dd con2}.
\end{proof}

Constraint \eqref{eq:stability1-design} is a linear matrix inequality (LMI). The nonconvexity comes from the equality constraint \eqref{eq:dd con1}. We will propose an iterative approach to efficiently deal with this nonconvexity in the next section.

% \subsection{Stabilization of nonlinear NFL (In progress)}
% Formulate the stabilization conditions of NFL as an SOS condition and transform it to an SDP problem. Then apply data-driven method to represent it.
% \begin{theorem}
% Here is my theorem.
% \end{theorem}

% \begin{proof}
% Here is my proof.
% \end{proof}

%%%%%%%%%%%%%%%%%%%%%%%%%%%%%%%%%%%%%%%%%%%%%%%%%%%%%%%%%%%%%%%%%%%%%%%%%%%%%%%%
\section{Iterative Algorithm for Stable NFL Design}\label{sec:Algorithm}
In this section we propose an iterative algorithm to solve Problem \ref{prob:1}. We incorporate the data-driven stability conditions \eqref{eq:stability-design} as hard constraints into an NN training framework with a certain loss function. 
% For the first algorithm we make the stability conditions strict while we lift the stability constraints \eqref{eq:stability1-design} into the objective function for the second algorithm. 
Based on this algorithm, a fine-tuning algorithm is proposed for solving Problem \ref{prob:2}.

\subsection{Iterative NFL Design Algorithm with Stability Constraints}

% {\color{blue}I think Section V.A and V.B can be merged.}

The stable NFL can be designed by solving the following optimisation problem.
% \begin{subequations}
% \label{eq:linear-optfortrain}
% \begin{gather}
%     \min_{N,Q,L} \eta_1\mathcal{L}(N)-\eta_2\log \det (Q_1)\
% \label{eq:obj-design}\\
%     s.t.  \text{     LMI}(Q,L) \succ 0\\
%     f(N)Q=\bar{U}_{0,T}L \label{eq:cons-1-design}\\
%     \begin{bmatrix}
%         I & 0\\
%     \end{bmatrix}Q 
%     =
%     \bar{X}_{0,T}L \label{eq:cons-2-design}
% \end{gather}
% \end{subequations}
\begin{subequations}\label{eq:linear-optfortrain}
    \begin{align}
        \min_{N,Q,L}& \eta_1\mathcal{L}(N)-\eta_2\log \det (Q_1)\
\label{eq:obj-design}\\\text{s.t.}~&Q_1\in\mathbb{S}_{++}^{n_x},Q_2\in\mathbb{S}_{++}^{n_\phi},\\
    &\eqref{eq:stability1-design},\eqref{eq:dd con2},\\
    &f(N)Q=\bar{U}_{0,T}L,\label{eq:linear-optfortrain-eq}
    \end{align}
\end{subequations}where 
$
    f(N) = \tilde N
$
represents the loop transformation.
The first term of the objective function is to minimise the loss function $\mathcal{L}(N)$ for NN prediction. It should be noted that we slightly abuse notation here: the prediction loss function $\mathcal{L}(\cdot)$ is actually a function of NN weights $W$ instead of matrix $N$, which is constructed from $W$. The second term is to maximise the ROA of the system. $\eta_1$ and $\eta_2$ are trade-off weights. Unlike traditional training process, we also add constraints from Theorem \ref{theo:stability_DDC} for stability guarantee.

To address the nonconvexity of optimisation problem \eqref{eq:linear-optfortrain}, we propose an algorithm to solve the problem iteratively. Following the framework proposed by \cite{RN15}, we place the nonconvex loop transformation equality constraints \eqref{eq:linear-optfortrain-eq} into the objective function using the augmented Lagrangian method. The augmented loss function is:
\begin{equation}\label{eq:nn-obj-func}
    \begin{split}
\mathcal{L}_a^1(N,Q,L,Y)=\eta_1\mathcal{L}(N)-\eta_2\log \det (Q_1)+ \\ 
    \text{tr}(Y^T(f(N)Q-\bar{U}_{0,T}L)+
    \frac{\rho}{2}||f(N)Q-\bar{U}_{0,T}L||^2_F
    \end{split}
\end{equation}
where $Y \in \mathbb{R}^{(n_u+n_{\phi})\times (n_x+n_{\phi})}$ is the Lagrange multiplier, and $\rho$ is a regularization parameter. 
% \begin{algorithm}
% \caption{Iterative design algorithm with hard constraints}
% \label{alg:LinearNFL1}
% \begin{algorithmic}[1] 
% \Procedure{LinearNFL1}{$U_{0,T},
% X_{0,T}, X_{1,T}$} 
% \State $k=0$, initialization
% \While {$||f(N^{k})Q^{k}-\bar{U}_{0,T}L^{k}||^2_F > \sigma$} 
%     \State $N$-update: 
%     \Statex\hspace{40pt}$N^{k+1}=\text{arg min}_N\mathcal{L}^1_a(N,Q^k,L^k,Y^k)$
%     \State ($Q,L$)-update:
%     \Statex\hspace{40pt}$(Q,L)^{k+1}=\text{arg min}_{Q,L}\mathcal{L}^1_a(N^{k+1},Q,L,Y^k)$
%     \Statex\hspace{40pt}s.t. 
%     \eqref{eq:stability1-design},\eqref{eq:dd con2}
    
% % LMI($Q,L)\succ0$, 
% % $\begin{bmatrix}
% %         I & 0\\
% %     \end{bmatrix}Q 
% %     =
% %     \bar{X}_{0,T}L$
%     \State $Y^{k+1}=Y^k+\rho(f(N^{k+1})Q^{k+1}-\bar{U}_{0,T}L^{k+1})$
%     \State $k$-update: $k=k+1$
% \EndWhile
% \State \Return the optimal result ($N^{k,*},Q^{k,*},L^{k,*},Y^{k,*}$) 
% \EndProcedure
% \end{algorithmic}
% \end{algorithm} 
\begin{algorithm}
\caption{Iterative design algorithm with hard constraints}
\label{alg:LinearNFL1}
\begin{algorithmic}[1] 
\Procedure{LinearNFL1}{$U_{0,T},
X_{0,T}, X_{1,T}$} 
\State $k=0$, initialization
\While {$||f(N^{k})Q^{k}-\bar{U}_{0,T}L^{k}||^2_F > \sigma$} 
    \State $N^{k+1}=\text{arg min}_N\mathcal{L}^1_a(N,Q^k,L^k,Y^k)$
    \State $(Q,L)^{k+1}=\text{arg min}_{Q,L}\mathcal{L}^1_a(N^{k+1},Q,L,Y^k)$
    \Statex\hspace{40pt}s.t. 
    \eqref{eq:stability1-design},\eqref{eq:dd con2}
    
% LMI($Q,L)\succ0$, 
% $\begin{bmatrix}
%         I & 0\\
%     \end{bmatrix}Q 
%     =
%     \bar{X}_{0,T}L$
    \State $Y^{k+1}=Y^k+\rho(f(N^{k+1})Q^{k+1}-\bar{U}_{0,T}L^{k+1})$
    \State $k=k+1$
\EndWhile
\State \Return the optimal result ($N^{k,*},Q^{k,*},L^{k,*},Y^{k,*}$) 
\EndProcedure
\end{algorithmic}
\end{algorithm} 

The solving process is shown in Algorithm \ref{alg:LinearNFL1}. Here $k$ is used to denote the number of iteration and $\sigma$ is used as a convergence criterion. In step 4, we use an IL framework to train the NN controller and update $N$. Then, we guarantee the stability of system by solving an SDP in step 5.

\subsection{Fine-Tuning Framework for Existing Unstable Neural Network Controller}
Another problem of interest is tuning an existing NN controller to render the closed-loop system stable. Due to the lack of stability guarantees for the traditional NN training process, the NFL may not be stable. Herein, we propose a verification and adaptation framework for existing NN controllers. First of all, stability of the NFL is verified by the data-driven stability verification algorithm proposed by our previous work \cite{wang2023model}. If the verification SDP is feasible, we conclude that the NFL is already stable. If we want to guarantee stability using this method, the NN controller will need to be fine-tuned, i.e. minimally adapted, to guarantee local stability.
The objective function for fine-tuning is shown as:
\begin{gather}\mathcal{L}^2_a(N_f,Q,L,Y)=\eta_3||N_f||^2_F-\eta_2\log \det (Q_1)+\notag \\ 
    \text{tr}(Y^T(f(\bar{N})Q-\bar{U}_{0,T}L)+
    \frac{\rho}{2}||f(\bar{N})Q-\bar{U}_{0,T}L||^2_F,
\end{gather}
where $N_f$ represents the amount that the NN should be fine-tuned. $\bar{N}=N+N_f$ contains weights and biases of the obtained NN controller after fine-tuning. This equation also guarantees that the original structure of the NN will not be changed. The first term in the objective function reflects the amount of tuning required, while the other terms are the same as \eqref{eq:nn-obj-func}. The fine-tuning algorithm is shown as Algorithm \ref{alg:Fine-tuning}.
\begin{algorithm}
\caption{Fine-tuning algorithm for an existing NN controller}
\label{alg:Fine-tuning}
\begin{algorithmic}[1] 
\Procedure{Fine-Tuning}{$N,U_{0,T}, X_{0,T}, X_{1,T}$} 
\State Data-driven stability verification \cite{wang2023model}
\If{verification SDP is feasible}
        \State \Return $\bar{N}=N, N_f=0$
    \Else
        \State $k=0$, initialization
\While {$||f(\bar{N}^{k})Q^{k}-\bar{U}_{0,T}L^{k}||^2_F > \sigma$} 
    \State $i=0$: 
    
    \While {$||N_f||^2_F > \sigma'$} 
    % \Statex\hspace{40pt}
    \State $\bar{N}^{k}_{i+1}=$
    \Statex\hspace{60pt}$\bar{N}^{k}_{i}+\text{arg min}_{N_f} \mathcal{L}_a^2(N_f,Q^k,L^k, Y^k)$
    \State Update $f(\bar{N}^k_{i+1})$ by \eqref{eq:f(N)compute} and \eqref{eq:Ccompute}
    \State $i=i+1$
    \EndWhile
    \State $\bar{N}^{k+1}=\bar{N}^k_i$
    \State $(Q,L)^{k+1}=$
    \Statex\hspace{80pt}$\text{arg min}_{Q,L}\mathcal{L}^2_a(\bar{N}^{k+1},Q,L,Y^k)$
    \Statex\hspace{45pt}s.t.  \eqref{eq:stability1-design},\eqref{eq:dd con2}
    \State $Y^{k+1}=Y^k+\rho(f(\bar{N}^{k+1})Q^{k+1}-\bar{U}_{0,T}L^{k+1})$
    \State $k=k+1$
\EndWhile
\State \Return the result ($\bar{N}^{k,*},Q^{k,*},L^{k,*},Y^{k,*}$) 
    \EndIf
\EndProcedure
\end{algorithmic}
\end{algorithm} 
% \begin{algorithm}
% \caption{Fine-tuning algorithm for an existing NN controller}
% \label{alg:Fine-tuning}
% \begin{algorithmic}[1] 
% \Procedure{Fine-Tuning}{$N,U_{0,T}, X_{0,T}, X_{1,T}$} 
% \State Data-driven stability verification \cite{wang2023model}
% \If{verification SDP is feasible}
%         \State \Return $\bar{N}=N, N_f=0$
%     \Else
%         \State $k=0$, initialization
% \While {$||f(\bar{N}^{k})Q^{k}-\bar{U}_{0,T}L^{k}||^2_F > \sigma$} 
%     \State $N$-update, $i=0$: 
    
%     \While {$||N_f||^2_F > \sigma'$} 
%     % \Statex\hspace{40pt}
%     \State $\bar{N}^{k}_{i+1}=$
%     \Statex\hspace{60pt}$\bar{N}^{k}_{i}+\text{arg min}_{N_f} \mathcal{L}_a^3(N_f,Q^k,L^k, Y^k)$
%     \State $f(\bar{N}^k_{i+1})$-update:
%     \State Update $f(\bar{N}^k_{i+1})$ by \eqref{eq:f(N)compute} and \eqref{eq:Ccompute}
%     \State $i$-update: $i=i+1$
%     \EndWhile
%     \State $\bar{N}^{k+1}=\bar{N}^k_i$
%     \State ($Q,L$)-update:
%     \Statex\hspace{40pt}$(Q,L)^{k+1}=\text{arg min}_{Q,L}\mathcal{L}^3_a(\bar{N}^{k+1},Q,L,Y^k)$
%     \Statex\hspace{45pt}s.t.  \eqref{eq:stability1-design},\eqref{eq:dd con2}
%     \State $Y^{k+1}=Y^k+\rho(f(\bar{N}^{k+1})Q^{k+1}-\bar{U}_{0,T}L^{k+1})$
%     \State $k$-update: $k=k+1$
% \EndWhile
% \State \Return the result ($\bar{N}^{k,*},Q^{k,*},L^{k,*},Y^{k,*}$) 
%     \EndIf
% \EndProcedure
% \end{algorithmic}
% \end{algorithm} 

Algorithm \ref{alg:LinearNFL1} uses a gradient descent method in an IL framework to update $N$. For Algorithm \ref{alg:Fine-tuning}, there is no need to retrain the NN. However, the objective function to be solved is nonconvex in $\bar{N}$. So we propose an iterative algorithm as shown in step 9 to step 13 to solve this nonconvex problem.

% {\color{blue}I got confused about the following statements, let's discuss.}
In the $(i+1)^{th}$ iteration, to eliminate the nonconvexity in the optimisation problem of step 10, we calculate the stacked sector bounds $A_{\phi}$, $B_{\phi}$ and $C_1,  C_2, C_3$ and $C_4$ based on the value of $\bar{N}^k_i$ from last iteration. This approximation is based on the assumption that the fine-tuning amount $\{N_f\}_i^k$ is small. Then $f(\bar{N}^k_{i+1})$ becomes a linear term on variable $\{N_f\}_{i+1}^k$, and the optimisation problem becomes a QP. After the first iterative algorithm converges, we update $Q, L$ in the same way as in Algorithm \ref{alg:LinearNFL1}.

% It should be noted that these quantities are all calculated from $\bar{N}_i^k=\bar{N}_{i-1}^k+\{N_f\}_{i-1}^k$, when the fine-tuning amount $\{N_f\}_i^k$ is small, we can assume these quantities remain unchanged. Therefore, in the $(i+1)^{th}$ iteration, we regard $f(\bar{N}^k_i)$ as an known parameter and use it to solve optimization problem. Then we update $f(\bar{N}^k_{i+1})$ based on new NN parameter $\bar{N}^k_{i+1}$. After the first iterative algorithm converges, we update $Q, L$.

%%%%%%%%%%%%%%%%%%%%%%%%%%%%%%%%%%%%%%%%%%%%%%%%%%%%%%%%%%%%%%%%%%%%%%%%%%%%%%%%
\section{Numerical Examples}\label{sec:case}
\subsection{Vehicle Controller Design}
To demonstrate the effectiveness of our method, we apply Algorithm \ref{alg:LinearNFL1} to the same numerical case as in \cite{RN15}, and compare the performance of our algorithm with the model-based NFL design approach in \cite{RN15}. The vehicle lateral dynamics have four states, $x=[e,\dot{e},e_{\theta},\dot{e}_{\theta}]^\top$, where $e$ and $e_{\theta}$ represents the perpendicular distance to the lane edge, and the angle between the tangent to the straight section of the road and the projection of the longitudinal axis of vehicle, respectively. $u\in\mathbb{R}$ is the steering angle of the front wheel.

For parameters, we use the constraint set $X=[-2,2]\times[-5,5]\times[-1,1]\times[-5,5]$. The NN controller for this system has two layers, each with 10 neurons. The expert demonstration data for IL training comes from an MPC law. Other parameters can be found in \cite{RN16}. To compare the proposed data-driven control (DDC) method with traditional model-based control (MBC) method, we use exactly the same data set for NN training. We choose $\rho=1000, \sigma=0.005$ and the maximum iteration number as 20 for the proposed algorithm. 

To demonstrate the effectiveness of the DDC method, we carry out several simulations with different weight combinations. After 10 experiments, the average prediction loss of NN trained with ($\eta_1=100,\eta_2=100$) is 0.128 while that with ($\eta_1=1000,\eta_2=100$) is 0.071. However, the ROA of the former is significantly larger than that of the latter. 

\begin{figure}[H]
  \centering
  \includegraphics[width=1\linewidth]{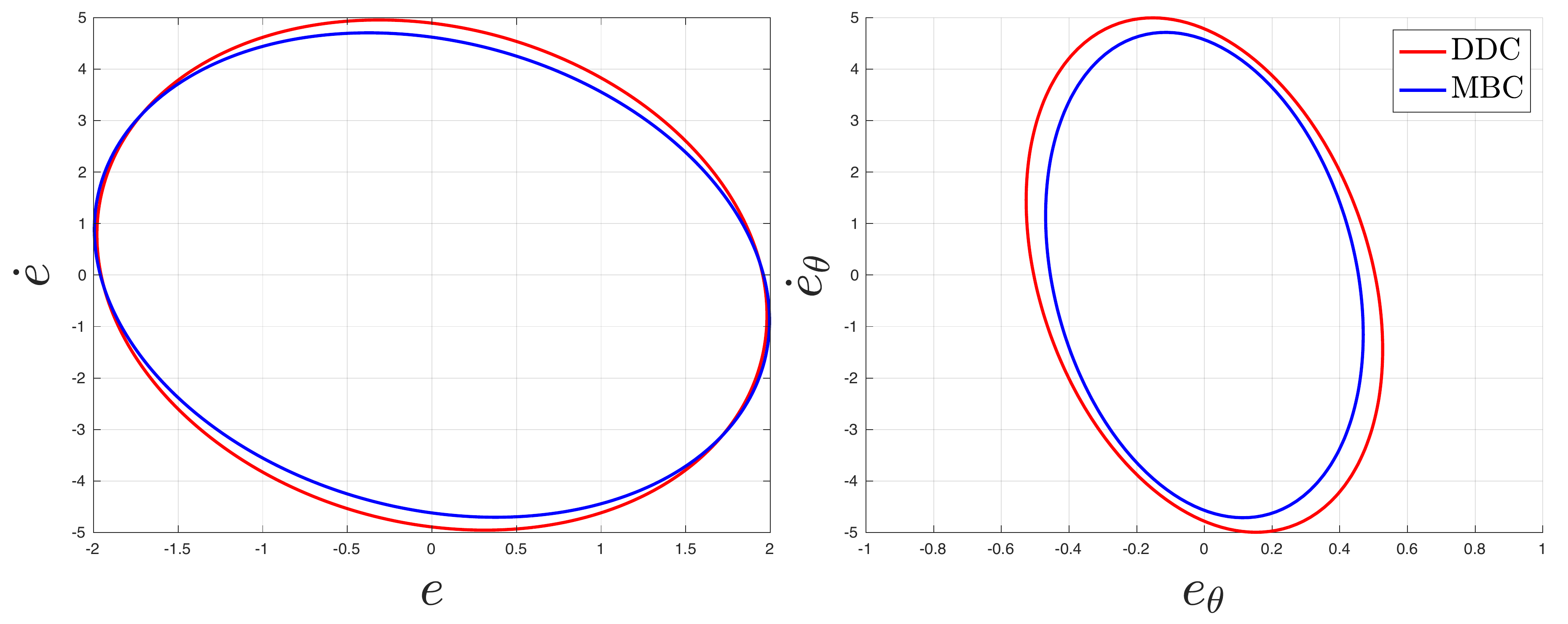}
  \caption{ROAs of NNs trained by DDC and MBC}\label{fig:roa}
\end{figure}

We then compare the proposed DDC approach with the MBC approach under the same settings: ($\eta_1=100,\eta_2=100$). The average prediction losses after 10 iterations is 0.128 for DDC, and 0.150 for MBC. Also, the ROAs under DDC are slightly bigger than that under MBC in most of the cases (7 times out of 10 experiments). ROAs under two approaches in one experiment are shown in Figure \ref{fig:roa}.
In this case, the NN trained by our proposed DDC approach shows better performance in terms of both prediction accuracy and the size of the ROA. This may mean that the controller design using data-driven methods can outperform that based on specific models. One intuitive explanation is that the traditional model-based method which contains system identification and control design can be
abstracted as a nested bilevel optimisation problem which usually leads to suboptimality. On the contrary, DDC is more flexible, making it more likely to find the ``optimal" solution for the control design. Readers interested in this are referred to relevant discussions in \cite{RN12, RN267}. 

% \begin{table}
% \caption{An Example of a Table}
% \label{table_example}
% \begin{center}
% \begin{tabular}{|c||c|}
% \hline
% One & Two\\
% \hline
% Three & Four\\
% \hline
% \end{tabular}
% \end{center}
% \end{table}
\subsection{Existing NN Controller Fine-Tuning}
We now consider the same task but in a different scenario: we assume that there is already a trained controller for the vehicle. As there is no stability guarantee in a traditional NN training algorithm, we verify its stability first. If it is proved to be unstable, we modify the NN to make the new NFL stable. As we have no access to the training data and hope its structure and other performances remain unchanged, the Algorithm \ref{alg:Fine-tuning} is applied to fine-tune it.

After fine-tuning the NN controller, we simulate the state variation from the same random initial point under input signals of two controllers. The sampling time is $0.02$ s and the result can be shown as following figure.
\begin{figure}[H]
  \centering
  \includegraphics[width=1\linewidth]{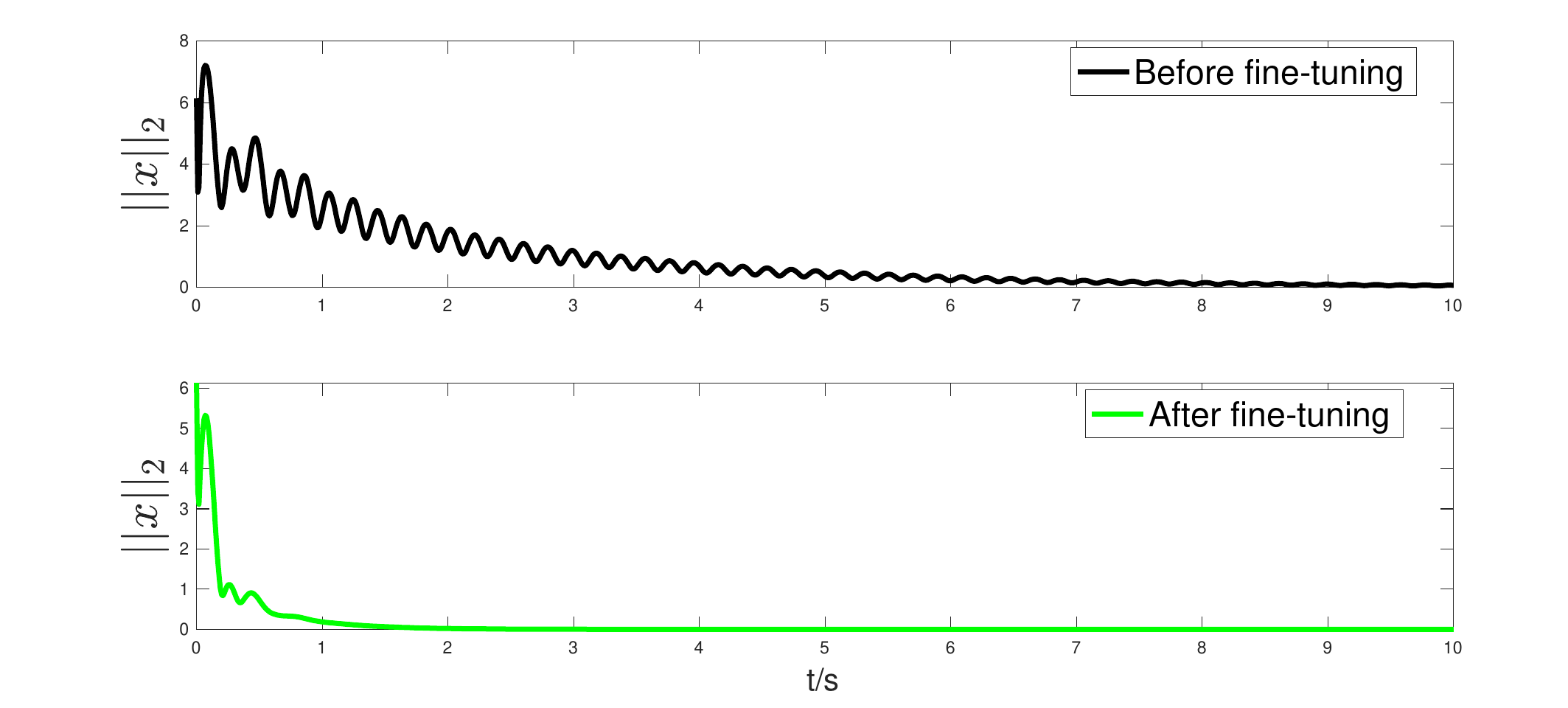}
  \caption{Simulation under the control of existing NN  and fine-tuned NN}\label{fig:simulation}
\end{figure}
Here we use $||x||_2$ as the value of vertical axis to reflect the variation of the state. Obviously, the NN controller before fine-tuning fails to stabilise the system in a long period, and it cannot be verified as stable by our verification method \cite{wang2023model}. After 5 iterations, the proposed algorithm converges and the fine-tuned controller is guaranteed to result in a stable NFL. The fine-tuning process is more efficient compared to training a new NN controller in terms of the running time. The total time for the former is 0.692 s and the average time for NN training is more than 30 mins. All simulations are performed on a laptop with Apple M2 chip. The MOSEK solver is used to solve the SDP and QP here. 
\section{Conclusions And Future Works}\label{sec:conclusion}
This paper proposes a data-driven method to design or fine-tune NN controllers for unknown linear systems. The design problem can be formulated as an optimisation problem and then solved by an iterative algorithm. Compared to traditional NN controller design methods, our algorithm solves an SDP that includes stability constraints during training and thus provides stability guarantees. The second algorithm we propose can provide stability guarantee for existing NN controllers by fine-tuning their parameters without changing their structure. Compared to retraining a controller for the system, this method significantly improves efficiency. We use a vehicle lateral control example to verify the proposed methods and compare it with model-based approaches. In our example, the controller designed directly shows better performance in terms of accuracy and the size of ROA compared to that designed based on model.

There are several promising directions for future works. First of all, we aim to extend our data-driven design method to NFLs involving nonlinear plants. Additionally, our method should be further developed to design more robust controllers to against disturbance and noise. Previous research has addressed this in LQR design problem by introducing different regularisers to achieve various data-driven parameterisations of the LQR \cite{RN255}. Finally, more theoretical works should be done to provide rigorous proof on the superiority of DDC for NN controllers design.

%%%%%%%%%%%%%%%%%%%%%%%%%%%%%%%%%%%%%%%%%%%%%%%%%%%%%%%%%%%%%%%%%%%%%%%%%%%%%%%%
% \section{ACKNOWLEDGMENTS}

% The authors gratefully acknowledge the contribution of National Research Organization and reviewers' comments.

%%%%%%%%%%%%%%%%%%%%%%%%%%%%%%%%%%%%%%%%%%%%%%%%%%%%%%%%%%%%%%%%%%%%%%%%%%%%%%%%

\bibliographystyle{IEEEtran} 
\bibliography{main} 

\begin{thebibliography}{10}
\providecommand{\url}[1]{#1}
\csname url@samestyle\endcsname
\providecommand{\newblock}{\relax}
\providecommand{\bibinfo}[2]{#2}
\providecommand{\BIBentrySTDinterwordspacing}{\spaceskip=0pt\relax}
\providecommand{\BIBentryALTinterwordstretchfactor}{4}
\providecommand{\BIBentryALTinterwordspacing}{\spaceskip=\fontdimen2\font plus
\BIBentryALTinterwordstretchfactor\fontdimen3\font minus \fontdimen4\font\relax}
\providecommand{\BIBforeignlanguage}[2]{{%
\expandafter\ifx\csname l@#1\endcsname\relax
\typeout{** WARNING: IEEEtran.bst: No hyphenation pattern has been}%
\typeout{** loaded for the language `#1'. Using the pattern for}%
\typeout{** the default language instead.}%
\else
\language=\csname l@#1\endcsname
\fi
#2}}
\providecommand{\BIBdecl}{\relax}
\BIBdecl

\bibitem{RN102}
A.~G. Barto, R.~S. Sutton, and C.~W. Anderson, ``Neuronlike adaptive elements that can solve difficult learning control problems,'' \emph{IEEE transactions on systems, man, and cybernetics}, no.~5, pp. 834--846, 1983.

\bibitem{RN101}
C.~W. Anderson, ``Learning to control an inverted pendulum with connectionist networks,'' in \emph{1988 American Control Conference}, 1988, Conference Proceedings, pp. 2294--2298.

\bibitem{RN72}
W.~Li and J.~J.~E. Slotine, ``Neural network control of unknown nonlinear systems,'' in \emph{1989 American Control Conference}, 1989, Conference Proceedings, pp. 1136--1141.

\bibitem{RN75}
Y.~Chen, Y.~Shi, and B.~Zhang, ``Optimal control via neural networks: A convex approach,'' \emph{arXiv preprint arXiv:1805.11835}, 2018.

\bibitem{RN76}
F.~Fabiani and P.~J. Goulart, ``Reliably-stabilizing piecewise-affine neural network controllers,'' \emph{IEEE Transactions on Automatic Control}, vol.~68, no.~9, pp. 5201--5215, 2023.

\bibitem{RN13}
E.~Kaufmann, L.~Bauersfeld, A.~Loquercio, M.~Müller, V.~Koltun, and D.~Scaramuzza, ``Champion-level drone racing using deep reinforcement learning,'' \emph{Nature}, vol. 620, no. 7976, pp. 982--987, 2023.

\bibitem{RN23}
S.~Gowal, K.~D. Dvijotham, R.~Stanforth, R.~Bunel, C.~Qin, J.~Uesato, R.~Arandjelović, T.~Mann, and P.~Kohli, ``On the effectiveness of interval bound propagation for training verifiably robust models,'' \emph{arXiv preprint}, 2019.

\bibitem{RN40}
M.~Fazlyab, M.~Morari, and G.~J. Pappas, ``Safety verification and robustness analysis of neural networks via quadratic constraints and semidefinite programming,'' \emph{IEEE Transactions on Automatic Control}, vol.~67, no.~1, pp. 1--15, 2022.

\bibitem{RN70}
M.~Newton and A.~Papachristodoulou, ``Neural network verification using polynomial optimisation,'' in \emph{2021 60th IEEE Conference on Decision and Control (CDC)}, 2021, Conference Proceedings, pp. 5092--5097.

\bibitem{RN39}
M.~Fazlyab, A.~Robey, H.~Hassani, M.~Morari, and G.~Pappas, ``Efficient and accurate estimation of lipschitz constants for deep neural networks,'' \emph{Advances in Neural Information Processing Systems}, vol.~32, 2019.

\bibitem{RN27}
P.~Pauli, A.~Koch, J.~Berberich, P.~Kohler, and F.~Allgöwer, ``Training robust neural networks using lipschitz bounds,'' \emph{IEEE Control Systems Letters}, vol.~6, pp. 121--126, 2022.

\bibitem{RN15}
H.~Yin, P.~Seiler, M.~Jin, and M.~Arcak, ``Imitation learning with stability and safety guarantees,'' \emph{IEEE Control Systems Letters}, vol.~6, pp. 409--414, 2022.

\bibitem{RN54}
F.~Gu, H.~Yin, L.~E. Ghaoui, M.~Arcak, P.~Seiler, and M.~Jin, ``Recurrent neural network controllers synthesis with stability guarantees for partially observed systems,'' \emph{Proceedings of the AAAI Conference on Artificial Intelligence}, vol.~36, no.~5, pp. 5385--5394, 2022.

\bibitem{RN28}
N.~Junnarkar, H.~Yin, F.~Gu, M.~Arcak, and P.~Seiler, ``Synthesis of stabilizing recurrent equilibrium network controllers,'' in \emph{2022 IEEE 61st Conference on Decision and Control (CDC)}, 2022, Conference Proceedings, pp. 7449--7454.

\bibitem{RN12}
R.~Sepulchre, ``Data-driven control: Part one of two [about this issue],'' \emph{IEEE Control Systems Magazine}, vol.~43, no.~5, pp. 4--7, 2023.

\bibitem{RN10}
C.~D. Persis and P.~Tesi, ``Formulas for data-driven control: Stabilization, optimality, and robustness,'' \emph{IEEE Transactions on Automatic Control}, vol.~65, no.~3, pp. 909--924, 2020.

\bibitem{RN71}
M.~Guo, C.~D. Persis, and P.~Tesi, ``Data-driven stabilization of nonlinear polynomial systems with noisy data,'' \emph{IEEE Transactions on Automatic Control}, vol.~67, no.~8, pp. 4210--4217, 2022.

\bibitem{RN87}
S.~Prajna, A.~Papachristodoulou, and W.~Fen, ``Nonlinear control synthesis by sum of squares optimization: a lyapunov-based approach,'' in \emph{2004 5th Asian Control Conference (IEEE Cat. No.04EX904)}, vol.~1, 2004, Conference Proceedings, pp. 157--165 Vol.1.

\bibitem{wang2023model}
H.~Wang, Z.~Xiong, L.~Zhao, and A.~Papachristodoulou, ``Model-free verification for neural network controlled systems,'' \emph{arXiv preprint arXiv:2312.08293}, 2023.

\bibitem{RN16}
H.~Yin, P.~Seiler, and M.~Arcak, ``Stability analysis using quadratic constraints for systems with neural network controllers,'' \emph{IEEE Transactions on Automatic Control}, vol.~67, no.~4, pp. 1980--1987, 2022.

\bibitem{RN267}
F.~Dörfler, J.~Coulson, and I.~Markovsky, ``Bridging direct and indirect data-driven control formulations via regularizations and relaxations,'' \emph{IEEE Transactions on Automatic Control}, vol.~68, no.~2, pp. 883--897, 2023.

\bibitem{RN255}
F.~Dörfler, P.~Tesi, and C.~De~Persis, ``On the role of regularization in direct data-driven lqr control,'' in \emph{2022 IEEE 61st Conference on Decision and Control (CDC)}.\hskip 1em plus 0.5em minus 0.4em\relax IEEE, 2022, Conference Proceedings, pp. 1091--1098.

\end{thebibliography}

\end{document}